\documentclass[12pt, final, 3p, times]{Qarticle}

\usepackage{epsfig}
\usepackage{amssymb}
\usepackage{float}
\usepackage{amsfonts}
\usepackage{mathrsfs}
\usepackage{amsthm}
\usepackage{amsmath}
\usepackage{amscd}
\usepackage{lineno}

\usepackage[colorlinks,citecolor=blue,linkcolor=blue,urlcolor=blue]{hyperref}
\usepackage{graphicx}
\usepackage{subfigure}

\usepackage{stmaryrd}
\usepackage{geometry}
\usepackage{fancyhdr}
\theoremstyle{plain}

\theoremstyle{definition}

\numberwithin{equation}{section}
\selectfont
\allowdisplaybreaks[4]

\begin{document}

\title{On an algorithm  for two-term spectral asymptotic formulas}
\fancyhead[C]{\footnotesize\textit{G. Q. Liu / On an algorithm  for two-term spectral asymptotic formulas}}

\begin{frontmatter}

\author[1]{Genqian Liu\corref{c}}
\ead{liuswgq@163.com}
 \cortext[c]{Corresponding author.}
 \address[1]{School of Mathematics and Statistics, Beijing Institute of Technology, Beijing 100081,  P. R. China}

\date{}

\begin{abstract} 
In the book [Yu. Safarov and D. Vassiliev, The asymptotic distribution of eigenvalues of partial differential operators, Amer. Math. Soc., Providence, RI, 1997], a key and central ``algorithm'' was established, by which the coefficients of two-term asymptotic expansions  of the  eigenvalue  counting functions can be explicitly calculated for many  partial differential operators under an additional geometric assumption. In this paper, we give a counter-example to this ``algorithm'' by discussing the case of elastic eigenvalues. This implies that the most conclusions in the above book written by Yu. Safarov and D. Vassiliev are fundamentally wrong because they are based on the erroneous ``algorithm''.
\end{abstract}


\end{frontmatter}

\section{Introduction}

Spectral asymptotics for partial differential operators have been the subject  of extensive research for over a century. It has attracted the attention of many  outstanding mathematicians and physicists (see, lines 2\! -- 4 on p.$\,$xi in \cite{SaVa-97}).
In order to study Weyl's conjecture concerning the existence of a
second asymptotic term  for the Laplacian, the Lam\'{e} operator,  and higher order partial differential operators,  Yu. Safarov and D. Vassiliev in \cite{SaVa-97} gave a strategy of  ``algorithm'' by which the coefficients of two-term asymptotic expansions can be explicit calculated under an  additional geometric assumption. This kind of ``algorithm'' as well as its proof is the most key and core part in \cite{SaVa-97}.

In two-term asymptotic formula of Weyl's conjecture for the Lam\'{e} operator,
 the authors of \cite{SaVa-97},  on p.$\,$236 of \S6.3, wrote \textcolor{blue}{\it ``Let $M$ be a region in $\mathbb{R}^n$, $y$ Cartesian coordinates
in $\mathbb{R}^n$, and $\mathbf{v}$ an $n$-component vector-function. We consider the spectral problem
for the system of equations}
\textcolor{blue}{\it \begin{eqnarray} \label{(A.1)}
  -c_t^2 \Delta \mathbf{v} -(c_l^2-c_t^2) \,\operatorname{grad}\; \operatorname{div}\;\mathbf{v} =\lambda^2 \mathbf{v}
 \end{eqnarray}
subject to the Dirichlet boundary conditions $\mathbf{v}\big|_{\partial M}=0$ (fixed boundary) or the
conditions of free boundary. The latter are the variational boundary conditions
generated by the quadratic functional
\begin{eqnarray} \label{(A.2)}
\mathcal{E} (\mathbf{v}):= \int_M \Big((c_l^2 -2 c_t^2) |\mbox{div}\;\mathbf{v}|^2 + \frac{c_t^2}{2} \sum_{i,j} |\partial_{y_j} \mathsf{v}_i +\partial_{y_i} \mathsf{v}_j|^2\Big) dy.  \end{eqnarray}}

\textcolor{blue}{\it The system (\ref{(A.1)}) describes the vibrations of an isotropic elastic body, see
[LanLif,\\ Sect.22]. Here $\lambda=\omega$ is the vibration frequency, and the constants $c_l$,
$c_t$ are the velocities of longitudinal and transverse waves, respectively. They are
assumed to satisfy the inequality $c_l/c_t >\sqrt{2}$. Further we denote $\alpha=c_t^2 c_l^{-2}$.\;\,............
 }

\textcolor{blue}{\it Weyl's formula for two-dimensional elasticity has the form
\begin{eqnarray}  \label{(A.3)}  N(\lambda) = \frac{(c_l^{-2} +c_t^{-2}) S}{4\pi} \lambda^2 +\frac{\beta L}{4\pi c_t} \lambda +o(\lambda), \;\;\; \;\;\;\, \lambda\to +\infty, \;\end{eqnarray} where $S$ is the surface area of $M$, $L$
 is the length of $\partial M$, and $\beta$ is a dimensionless
coefficient. The value of the coefficient $\beta$ for a fixed boundary is
\begin{eqnarray} \label{(A.4)}  \beta= -1 -\sqrt{\alpha} -\frac{4}{\pi} \int_{\sqrt{\alpha}}^1 \tan^{-1} \sqrt{(1-\alpha \xi^{-2})(\xi^{-2} -1) }\, d\xi. \end{eqnarray}
For free boundary
\begin{eqnarray}  \label{(A.5)} \beta = 4\gamma^{-1} -3 +\sqrt{\alpha} +\frac{4}{\pi} \int^1_{\sqrt{\alpha}} \tan^{-1} \frac{(2-\xi^{-2})^{2}}{4\sqrt{(1-\alpha \xi^{-2} )(\xi^{-2} -1)}}\,d\xi, \qquad \end{eqnarray}
where $0<\gamma<1$ is the root of the algebraic equation
\begin{eqnarray}  \label{(A.6)}   \gamma^6 -8\gamma^4 +8(3-2\alpha )\gamma^2 -16(1-\alpha)=0.\;''  \end{eqnarray}}

In this paper,  we will show that in  the two-term asymptotic expansion (\ref{(A.3)}), the constant $\beta$ in  (\ref{(A.5)})  is wrong because the corresponding result with such a constant  contradicts to the generalized  McKean--Singer's conclusion in \cite{MS-67}. This implies that the ``algorithm''  (i.e., the central theory) in \cite{SaVa-97} is  fundamentally wrong,  which leads to  that  the most  conclusions in \cite{SaVa-97} are  essentially
  wrong. Finally, we will also  give the  correct  two-term asymptotic expansions of the elastic heat traces as well as the corresponding eigenvalue counting functions in a smooth compact Riemannian manifold with  boundary.

\section{A counter-example to the algorithm method}

It is well-known (see, for example,  Chapter I of \cite{KGBB}) that for the  system (\ref{(A.1)}) of elastic eigenvalue equations, the free boundary condition (i.e., traction)  is $$2c_t^2\,\big(\mbox{Def}\; \mathbf{v}\big)\,\boldsymbol{\nu}+(c^2_l-2c_t^2) (\mbox{div}\; \mathbf{v})\,\boldsymbol{\nu}=0 \;\;\;\mbox{on}\;\;\,\partial M,$$
where $\mathbf{v}=(\mathsf{v}_1, \cdots, \mathsf{v}_n)\in [C^2(\bar M)]^n$,  \begin{eqnarray*}\mbox{Def}\; \mathbf{v} :=\frac{1}{2}\big[\partial_{y_j} \mathsf{v}_i + \partial_{y_i}\mathsf{v}_j\big]_{n\times n}=\frac{1}{2}\begin{bmatrix} \partial_{y_1} \mathsf{v}_1+\partial_{y_1} \mathsf{v}_1 & \partial_{y_2} \mathsf{v}_1 +\partial_{y_1} \mathsf{v}_2 &\cdots & \partial_{y_n} \mathsf{v}_1 +\partial_{y_1} \mathsf{v}_n\\
\partial_{y_1} \mathsf{v}_2+\partial_{y_2} \mathsf{v}_1 & \partial_{y_2} \mathsf{v}_2 + \partial_{y_2} \mathsf{v}_2&\cdots & \partial_{y_n} \mathsf{v}_2+ \partial_{y_2} \mathsf{v}_n\\
\cdots & \cdots & \cdots\\
  \partial_{y_1} \mathsf{v}_n + \partial_{y_n} \mathsf{v}_1& \partial_{y_2} \mathsf{v}_n + \partial_{y_n} \mathsf{v}_2& \cdots & \partial_{y_n} \mathsf{v}_n+\partial_{y_n} \mathsf{v}_n\end{bmatrix}, \end{eqnarray*}
and $\boldsymbol{\nu}=(\nu_1,\cdots, \nu_n)$ is the unit outer normal to $\partial M$. This free boundary condition can be directly obtained by Green's formula (see,  \S1 of Chapter III in \cite{KGBB}).

On p.$\,$236 of \S6.3 in \cite{SaVa-97}, the Lam\'{e} coefficients $c_t^2$ and $c_l^2-2c_t^2$ were restricted to $c_t^2>0$ and $c_l^2-2c_t^2>0$, in other words,  $0<\alpha=c_t^2c_l^{-2}<\frac{1}{2}$. This is
a very strict restriction for the Lam\'{e} coefficients, which excludes many useful physical and mathematical models. In \cite{CaFrLeVa-23},   which was a generalization of (\ref{(A.3)}) -- (\ref{(A.5)}) in a Riemannian manifold,  the range of the  Lam\'{e} coefficients was relaxed to  $c_t^2>0$ and $c_l^2-c_t^2>0$ (i.e., $0<\alpha<1$).  Obviously, the authors of \cite{SaVa-97} and \cite{CaFrLeVa-23} were not clear about the true range of the Lam\'{e} coefficients because they had wrongly thought that the ellipticity  of the corresponding boundary value problems might be broken if $c_t^2$, $c_l^2$ lie outside the range $c_t^2>0$, $c_l^2-c_t^2>0$  as remarked  in footnote 5 on page 3 of \cite{CaFrLeVa-23}  (Let us point out in advance that the true range of the Lam\'{e} coefficients is $c_t^2>0$ and $c_l^2>0$, which ensures the  Lam\'{e}  system $c_t^2 {\Delta} \mathbf{I}_n  +(c_l^2-c_t^2) \, \mbox{grad}\,\mbox{div}$ is strongly elliptic).

For the $n$-dimensional elastic Lam\'{e} system ($n\ge 2$), the Lam\'{e} coefficients can be further extended to the {\it strong convexity condition} $c_t^2>0$ and $n(c_l^2-2c_t^2) + 2 c_t^2>0$ (see, for example, p.$\,$169 of \cite{Isa-06} for $n$-dimensional case,  and (5.171) on p.213 of \cite{LiQin-13} or (1.11) on p.$\,$46 of \cite{KGBB} for three-dimensional case), in which the Lam\'{e} coefficients have physical explanation (see,  p.$\,$11 of \cite{LanLif}). This is also the main range which can be  studied by physics and mechanics experts by variational methods during an early stage. Thus in the two-dimensional case, one has $0<\alpha=c_t^2c_l^{-2}<1$. Note that for $n=2$, when $c_t^2>0$ and $c_l^2-c_t^2=0$, it can be directly verified that the corresponding elastic energy density in (\ref{(A.2)}) is pointwise non-negative in $M$.

In the discussions about the systems of elastostatics (i.e., steady-state case), one of the fundamental questions is how to give assumptions on the elastic tensor, such that these systems can be mathematically and effectively analyzed. Therefore,  in mathematics  (particularly, in the field of partial differential equations), it is always assumed that for any dimensions $n\ge 2$, the Lam\'{e} coefficients $c_t^2$ and $c_l^2-2c_t^2$ satisfy (see, for example,   \cite{CiMa},  (10.4) on p.$\,$297 of \cite{McL}, Chapters 5--6 of \cite{Hah-98},  (5.150) of \cite{LiQin-13},  \cite{MaHu} or \cite{MHNZ}):
\begin{eqnarray}  \label{(A.7)}  c^2_t>0 \;\;\mbox{and}\;\; (c_l^2-2c_t^2) +2c_t^2=c_l^2>0.\end{eqnarray}
The above condition (\ref{(A.7)}) is the most correct (and precise) range for the Lam\'{e} coefficients in studying the Lam\'{e} system of elastic equations.
In fact, from $a_{ijkl}= (c_l^2-2c_t^2)  \delta_{ij}\delta_{kl} + c_t^2 (\delta_{ik} \delta_{jl}+\delta_{il}\delta_{jk})$, one can
 verify that the  condition  (\ref{(A.7)}) of the Lam\'{e} coefficients $c_t^2$ and $c_l^2-2c_t^2$ is just equivalent to the following strong ellipticity condition (see p.$\,$208 of \cite{LiQin-13} or p.$\,$241 of \cite{MaHu}):

\vskip 0.19 true cm

 \noindent{\bf Proposition 2.1.} { \it If there exists a constant $\Xi>0$ such that
\begin{eqnarray*} \sum_{i,j,k,l=1}^n a_{ijkl} \xi_i\xi_k \eta_j \eta_l\ge \Xi \,|\xi|^2 |\eta|^2, \;\;\;
\forall \xi, \eta\in \mathbb{R}^n,\end{eqnarray*}
 then the fourth-order tensor $\mathbf{A}=(a_{ijkl})$ satisfies the strong ellipticity condition}.

\vskip 0.29 true cm 
If $c_t^2>0$ and $c_l^2>0$, then the isotropic strongly elliptic tensor $a_{ijkl}=(c_l^2-2c_t^2)  \delta_{ij}\delta_{kl} + c_t^2 (\delta_{ik} \delta_{jl}+\delta_{il}\delta_{jk})$ has clearly physical meaning (see,  p.$\,$242 of \cite{MaHu}):  
there are two speeds of propagation of plane progressive elastic waves given by 
\begin{eqnarray*} c_1=\sqrt{\frac{c_l^2} {\rho}} \;\;\, \mbox{and}\;\;\, c_2=\sqrt{\frac{c_t^2}{\rho}}\end{eqnarray*}
($c^2$ is a double eigenvalue).

Let us point out that for the elastic operator $${P}\mathbf{v}:=-c_t^2\Delta \mathbf{v} - (c_l^2-c_t^2) \,\mbox{grad}\; \mbox{div}\; \mathbf{v}, \;\;\;\,\;\; \mathbf{v}=(\mathsf{v}_1,\cdots, \mathsf{v}_n)$$  with   $c_t^2>0$ and $c_l^2>0$,
 the above definition of the strong ellipticity condition  is the same as given  by the classical ellipticity condition in (2.5) on p.$\,$44 of \cite{ADN-64} for the differential operator ${P}$ (see also,  p.$\,$193, (6.5) of \cite{McL}). Clearly, if $c_t^2>0$ and $c_l^2> 0$, then the  Lam\'{e} system of  elastic equations  $$-c_t^2\Delta \mathbf{v} - (c_l^2-c_t^2) \,\mbox{grad}\; \mbox{div}\; \mathbf{v}=\mathbf{f} $$ is strong ellipticity,
and the Dirichlet and free boundary conditions both satisfy the complementing boundary conditions (see, p.$\,$42 -- 44 of \cite{ADN-64}), where $\mathbf{f} =-\rho \mathbf{b} -\mbox{div}\, \boldsymbol{\sigma}$, $\mathbf{b}$ is the body force, $\boldsymbol{\sigma}$  is the Cauchy stress in the linearizing configuration  (see, Chapter 6, on p.$\,$316 in \cite{MaHu}).
 In other words, the ellipticity of the corresponding (Dirichlet and free) boundary conditions can make these two types boundary value   problems of the Lam\'{e} system  to be well-posed if the Lam\'{e} coefficients $c_t^2$ and $c_l^2- 2c_t^2$ lie
inside the range $c_t^2 >0, \, c_l^2> 0$ (see, p.$\,$331 of \cite{MaHu}).

  In eigenvalue problems, the most correct (and precise) range of the Lam\'{e} coefficients $c_t^2$ and $c_l^2-2c_t^2$ is also $c_t^2>0, \,c_l^2>0$, which implies  $0<\alpha=c_t^2c_l^{-2}<+\infty$.

\vskip 0.22 true cm

On the other hand,  the  {\it generalized Ahlfors-Laplacian} (see,
\cite{BGOP} or \cite{PiOr-96}),  which  stems from the study of quasiconformal geometry,  is the operator defined on $1$-form of a Riemannian manifold $M$:
 \begin{eqnarray} \label{023.9.21-2} A:=a \, d\delta  + b\,\delta d  - 2b \, \mbox{Ric},\end{eqnarray}
  where $a$ and $b$  are positive constants,  $\mbox{Ric\,}$ denotes the Ricci action on $1$-forms,
   $d$ is the exterior differential operator, and $\delta$ the adjoint operator of  $d$. 
  It has been shown in \cite{Liu-23} that  the Lam\'{e} operator  $P$, which is defined on vector field, can be equivalently  written as the {\it generalized Ahlfors-Laplacian} $A$ with $a:=c_l^2$ and $b:=c_t^2$, which is defined on $1$-form. This clearly indicates that the range of the Lam\'{e} 
   coefficients is  $c_t^2>0$ and $c_l^2>0$

\vskip 0.28 true cm 

For convenience, let ${H}_\partial^2$ denote the set of $\mathbf{v}\in [H^2(M)]^n$ with homogeneous boundary conditions of displacement or traction (or both) imposed; that is,
\begin{eqnarray*} \mathbf{v}=0 \;\, \mbox{on}\;\, \partial_{{D}}\;\;\, \mbox{and}\;\; (\mathbf{A}\cdot \nabla \mathbf{v})\cdot \boldsymbol{\nu}:= \Big( 2c_t^2
 (\mbox{Def}\; \mathbf{v}) \,\mathbf{I}_n +(c_l^2-2c_t^2) \, \mbox{div}\; \mathbf{v} \Big) \,\boldsymbol{\nu}=0\;\;\mbox{on}\;\, \partial_{{T}},\end{eqnarray*}
where $\mathbf{A}=(a_{ijkl})= (c_l^2-2c_t^2)  \delta_{ij}\delta_{kl} + c_t^2 (\delta_{ik} \delta_{jl}+\delta_{il}\delta_{jk})$ is the fourth-order tensor,  $\partial_{D}\subset \partial M$ and $\partial_T\subset \partial M$ are the subdomain of $\partial M$,  $\;\partial_D\cap \partial_T\ne \emptyset$,  $\;\partial_D \cup \partial_T=\partial M$,   and  one of $\partial_D$ or $\partial_T$ may be empty (cf. p.$\,$317 of \cite{MaHu}).  When $\partial_T=\emptyset $,  this is just the Dirichlet boundary condition,  while when   $\partial_D=\emptyset$, this exactly  becomes the free boundary condition.   Define  the linear operator
\begin{eqnarray*} {P}_\partial: H_\partial^2 \to L^2\;\;\;\;\mbox{by}\;\;\; {P}_\partial (\mathbf{v}) =\mbox{div} (\mathbf{A}\cdot \nabla\mathbf{v})= -c_t^2 \Delta \mathbf{v} -(c_l^2-c_t^2) \,\mbox{grad}\; \mbox{div}\;\mathbf{v} .\end{eqnarray*}
Denote by $$B(\mathbf{v},\mathbf{v})=\int_M \big((c_l^2 -2c_t^2) (\mbox{div}\; \mathbf{v})(\mbox{div}\; \mathbf{v}) +2c_t^2 \big\langle \mbox{Def}\; \mathbf{v}, \mbox{Def}\; \mathbf{v}\big\rangle \big)dy$$ the Dirichlet form corresponding to the elastic operator $P_\partial$ on space $H^1\times H^1$ (see, (3) on p.$\,$317 in \cite{MaHu}), where $\langle \mbox{Def}\; \mathbf{v}, \mbox{Def}\; \mathbf{v}\rangle =\frac{1}{4}\sum_{i,j} |\partial_{y_j} \mathsf{v}_i +\partial_{y_i} \mathsf{v}_j|^2$.

\vskip 0.1 true cm
The following proposition plays a key role for the Lam\'{e} operator $P_\partial$:

\vskip 0.18 true cm
\noindent {\bf  Proposition 2.2 (see p.$\,$318 of \cite{MaHu})}.\   {\it Let $P_\partial$ be defined above, and let $P_\partial$ be symmetric.  Then strong
ellipticity is equivalent to G{\aa}rding's inequality:
there are constants $\tilde{c}>0$ and $\tilde{d}>0$ such that for all $\mathbf{v}\in H^1$,
\begin{eqnarray} \label{23.12.10-4} B(\mathbf{v},\mathbf{v})\ge \tilde{c}\|\mathbf{v}\|_{H^1}^2 -\tilde{d} \|\mathbf{v}\|_{L^2}^2.\end{eqnarray}}

\vskip 0.42 true cm

 Now, we can verify that ${P}_\partial$
is symmetry; that is,
\begin{eqnarray*} \langle\langle {P}_\partial \mathbf{u} , \mathbf{v}\rangle\rangle_{L^2}= \langle\langle \mathbf{u}, P_\partial \mathbf{v}\rangle\rangle_{L^2}\;\;\, \mbox{for all}\;\, \mathbf{u}, \mathbf{v} \in H_\partial^2,
 \end{eqnarray*}
 where $\langle\langle \,\cdot  , \cdot \,\rangle\rangle_{L^2} $ is the $L^2$ inner product, and $ \langle \cdot,\cdot \rangle$ is the usual inner product in the sense of the Euclidean space.
In fact, for $\mathbf{u}, \mathbf{v}\in H_\partial^2$, we have
\begin{eqnarray*}&& \langle\langle P_\partial \mathbf{u}, \mathbf{v}\rangle\rangle_{L^2} =  \int_{M}\big \langle\!-(c_l^2-c_t^2)\, \mbox{grad}\; \mbox{div}\;\mathbf{u} -  c_t^2 \,\Delta \mathbf{u}, \mathbf{v}\big\rangle\,dy \\
&&  + \int_{\partial M}  \big\langle (c_l^2-2c_t^2) (\mbox{div}\;\mathbf{u}\big) \boldsymbol{\nu} + 2c_t^2 (\mbox{Def}\, \mathbf{u}\big) \boldsymbol{\nu} ,\mathbf{v}\big\rangle \,dS_y\\
&&=\int_{M}\big \langle-(c_l^2-2c_t^2)\, \mbox{grad}\; \mbox{div}\;\mathbf{u} + 2 c_t^2 \,\big(-\frac{1}{2}\Delta \mathbf{u} - \frac{1}{2} \,\mbox{grad}\;\mbox{div}\; \mathbf{u}\big), \mathbf{v}\big\rangle \,dy\\
&& + \int_{\partial M}  \big\langle (c_l^2-2c_t^2) (\mbox{div}\;\mathbf{u}\big) \boldsymbol{\nu} + 2c_t^2 (\mbox{Def}\, \mathbf{u}\big) \boldsymbol{\nu} ,\mathbf{v}\big\rangle dS_y
\\
&&= \int_M \big\langle -(c_l^2-2c_t^2) \,\mbox{grad}\;\mbox{div}\; \mathbf{u} + 2 c_t^2 \,\mbox{Def}^* \,\mbox{Def} \; \mathbf{u}, \mathbf{v}\big\rangle dy\\
&& + \int_{\partial M} \big\langle (c_l^2-2c_t^2) (\mbox{div}\; \mathbf{u}) \, \boldsymbol{\nu} + 2c_t^2 (\mbox{Def}\; \mathbf{u} ) \, \boldsymbol{\nu}, \mathbf{v} \big\rangle dS_y\\
&&= \int_M (c_l^2 -2c_t^2) (\mbox{div}\; \mathbf{u})(\mbox{div}\; \mathbf{v}) +2c_t^2 \big\langle \mbox{Def}\; \mathbf{u}, \mbox{Def}\; \mathbf{v}\big\rangle dy\\
&& =
 \int_{M}\big \langle\mathbf{u}, -(c_l^2-c_t^2)\, \mbox{grad}\; \mbox{div}\;\mathbf{v} -  c_t^2 \,\Delta \mathbf{v} \big\rangle \\
&& + \int_{\partial M}  \big\langle \mathbf{u}, (c_l^2-2c_t^2) (\mbox{div}\;\mathbf{v}\big) \boldsymbol{\nu} + 2c_t^2 (\mbox{Def}\, \mathbf{v}\big) \boldsymbol{\nu}\big\rangle dS_y
 =\langle\langle \mathbf{u}, P_\partial \mathbf{v}\rangle\rangle_{L^2}. \end{eqnarray*}
Here we have used the corresponding boundary conditions and the fact (see p.$\,$562 in \cite{Ta3}) that \begin{eqnarray*} \mbox{Def}^* \,\mbox{Def}\;\mathbf{u}:=-\mbox{div}\, \mbox{Def}\; \mathbf{u}=-\frac{1}{2} \Delta \mathbf{u} -\frac{1}{2} \mbox{grad}\; \mbox{div}\; \mathbf{u}. \end{eqnarray*}
  Moreover, on p.$\,$345-346 of \cite{MaHu}, Marsden and Hughes also proved that the equations of (linear elasticity) motion generate a quasi-contractive semigroup in
$H^1\times L^2$ (relative to some Hilbert space structure) under the assumption $c_t^2>0$ and $c_l^2>0$.

For $c_t^2>0$ and $c_l^2>0$, one can consider the following elastic eigenvalue problems:
\begin{eqnarray} \label{(A.9)}  \left\{ \begin{array}{ll} -c_t^2 \Delta \mathbf{v} -(c_l^2-c_t^2) \,\mbox{grad}\; \mbox{div}\;\mathbf{v} =\tau \mathbf{v} \;\;\; &\mbox{in}\;\; M,\\
 \mathbf{v}=0 \;\;\;&\mbox{on}\;\; \partial M \end{array} \right.
\end{eqnarray}
and
\begin{eqnarray} \label{(A.10)}  \left\{ \begin{array}{ll} -c_t^2 \Delta \mathbf{v} -(c_l^2-c_t^2) \,\mbox{grad}\; \mbox{div}\;\mathbf{v} =\tau \mathbf{v} \;\;\; &\mbox{in}\;\; M,\\
 c_t^2 \big[\partial_{y_j} \mathsf{v}_i + \partial_{y_i}\mathsf{v}_j\big]_{n\times n}\,\boldsymbol{\nu}+(c^2_l-2c_t^2) (\mbox{div}\; \mathbf{v})\,\boldsymbol{\nu}=0 \;\;\;&\mbox{on}\;\; \partial M. \end{array} \right.
\end{eqnarray}
Here, for the sake of convenience, in (\ref{(A.9)}) and (\ref{(A.10)}) we have written the eigenvalue as $\tau$ instead of $\lambda^2$.
Let us take 
\begin{eqnarray*}\interleave \mathbf{v}  \interleave = \left( \int_M \Big( \big( c_l^2 -2c_t^2) \,\big(\mbox{div}\, \mathbf{v}\big) \big(\mbox{div}\, \mathbf{v}\big)   + 2 c_t^2 \langle \mbox{Def}\, \mathbf{v}, \mbox{Def}\, \mathbf{v}\rangle  + \tilde{d} \, |\mathbf{v}|^2 \Big) \,dy \right)^{1/2},\end{eqnarray*}
$\tilde{d}$ is the corresponding (second term) constant in G{\aa}rding's inequality (\ref{23.12.10-4}).  Then, by G{\aa}rding's inequality we see that there exists constant $\tilde{d}_1>0$ such that 
\begin{eqnarray*} \sqrt{\tilde{c}}\;\| \mathbf{v} \|_{H^1} \le   \interleave \mathbf{v}  \interleave \le \tilde{d}_1\, \| \mathbf{v}\|_{H^1}, \;\;\, \quad \forall\, \mathbf{v}\in H_\partial^1 (M), \end{eqnarray*}
where $\tilde{c}$ is the (first term) constant in G{\aa}rding's inequality  (\ref{23.12.10-4}).  Clearly, $\interleave \cdot  \interleave$ is an equivalent normal in $H^1_\partial (M)$, so the $H_\partial^1(M)$ is a Banach space under the norm $\interleave \cdot  \interleave$. Put 
\begin{eqnarray*}  (\mathbf{u},\mathbf{v})_{\tau_0} =\int_{M}\Big( B(\mathbf{u},\mathbf{v}) \,+ \tau_0 \langle \mathbf{u}, \mathbf{v}\rangle \Big) \,dy,\end{eqnarray*} 
we see that  $H_\partial^1$ is a Hilbert space under $(\cdot, \cdot)_{\tau_0}$,  where $\tau_0=\tilde{d}$ (the second  constant in 
G{\aa}rding's inequality  (\ref{23.12.10-4})). 
It follows that for any $\mathbf{u}\in L^2(M)$, we see 
\begin{eqnarray*}  \bigg |\int_M \langle \mathbf{u}, \mathbf{v}\rangle   \,dy \bigg|  \le \| \mathbf{u}\|_{L^2(M)} \|\mathbf{v}\|_{L^2(M)}  \le \tilde{c}_1 \,\|\mathbf{u}\|_{L^2(M)} \interleave \!\mathbf{v}  \interleave,  \;\;\;\;\forall \,\mathbf{v}\in (H^1_\partial (M), \interleave \cdot \interleave) ,  \end{eqnarray*} where $\tilde{c}_1=1/\!\sqrt{\tilde{c}}$.
According to the Riesz theorem, there exists unique $\mathbf{w}\in (H_\partial^1(M), \interleave \cdot \interleave)$ such that 
\begin{eqnarray} \label{23.10.17-1}\int_M \langle\mathbf{u},\mathbf{v}\rangle \, dy =( \mathbf{w}, \mathbf{v})_{\tau_0}, \;\;\;\forall \, \mathbf{v}\in (H_\partial^1(M),\interleave \cdot \interleave).\end{eqnarray} 
Define $K_{\tau_0} : ( L^2(M) , \|\cdot\|_{L^2}) \to (H_\partial^1 (M), \interleave \cdot\interleave)$ by $\mathbf{w} :=K_{\tau_0} \mathbf{u}$. Then 
\begin{eqnarray*} \interleave K_{\tau_0} \mathbf{u}\interleave \le \tilde{c}_1 \|\mathbf{u}\|_{L^2(M)}, \;\;\, \forall \, \mathbf{u} \in L^2(M),\end{eqnarray*}
which implies that $K_{\tau_0}: \, \big( L^2(M),\|\cdot\|_{L^2}\big) \to \big( H_\partial^1 (M),\interleave\cdot\interleave) $ is a continuous linear operator. Denote by $ \iota$ the embedding operator form $(H_\partial^1(M),\interleave \cdot\interleave)$ to $(L^2(M), \|\cdot\|_{L^2})$. Rellich's theorem implies that $\iota$ is a compact operator. Thus the eigenvalue problem (\ref{(A.9)}) or (\ref{(A.10)})
 is equivalent to the following equations 
 \begin{eqnarray*} \left( I- (\tau +\tau_{0}) K_{\tau_0} \iota\right)
\mathbf{v} =\mathbf{0}, \quad \quad  \,\mathbf{v}\in (H^1_\partial (M), \interleave\cdot\interleave),\end{eqnarray*} 
where $I$ is the identity operator on $ (H^1_\partial (M), \interleave\cdot\interleave)$.
By applying the Riesz-Schauder theory and Hilbert-Schmidt theory, we know that $\sigma(K_{\tau_0} \iota)\setminus \{0\}$ is the set of real numbers and it contains  countable real numbers. 
Thus, there exist sequences $\{(\tau_k^{\begin{small}\mbox{Dir}\end{small}}, \mathbf{v}_k^{\mbox{Dir}})\}_{k=1}^\infty$ and $\{(\tau_k^{\mbox{free}}, \mathbf{v}_k^{\mbox{free}})\}_{k=1}^\infty$ of elastic eigen-pairs  satisfying (\ref{(A.9)}) and (\ref{(A.10)}), respectively.
It is east to verify that 
\begin{eqnarray*} \tau_k^{\mbox{Dir}\,/\,\mbox{free}}=&& \!\!\!\!\!\sup\limits_{E_{k-1}} \, \inf\limits_{\underset{\mathbf{v\ne 0}}{\mathbf{u}\in E_{k-1}^{\bot} } }
\frac{(\mathbf{v},\mathbf{v})_{\lambda_0}}{(K_{\lambda_0}\iota \,\mathbf{v}, \mathbf{v})_{\lambda_0}} -\tau_0
\\
=&&\!\!\!\!\!\sup\limits_{E_{k-1}} \, \inf\limits_{\underset{\mathbf{v}\ne 0}
{\mathbf{v} \in E_{k-1}^\bot}} \frac{\int_M B(\mathbf{v}, \mathbf{v}) \, dy }{\int_M |\mathbf{v}|^2 dy},\end{eqnarray*} 
where $E_{k-1}$ is any closed sub-space of $H_\partial^1 (M)$ of dimension $k-1$,  $\,\; k=1, 2, \cdots$ (For the Dirichlet eigenvalue problem  we should take $\partial_D=\partial M$, and for the free boundary eigenvalue problem we should take $\partial_T=\partial M$).  \  In particular, for any fixed positive  integer $k\ge 1$, the eigenvalues $\tau_k^{\mbox{Dir}}$ and $\tau_k^{\mbox{free}}$ are continuous functions in variable $\eta:=c_l^2-c_t^2$. For clarity, let us denote  by  $\tau_k^{\mbox{Dir}}(c_l^2-c_t^2)$ and $\tau_k^{\mbox{free}}(c_l^2-c_t^2)$ the eigenvalues  $\tau_k^{\mbox{Dir}}$ and $\tau_k^{\mbox{free}}$, respectively.

In \cite{SaVa-97}, the eigenvalue counting function $N(\lambda)$ for (\ref{(A.1)}) was defined by $N(\lambda):=\max\limits_{\substack{ k}} \{ k\big| \lambda_k<\lambda\}$. However, when $\lambda_k^2$ is replaced by $\tau_k$, we will use another
 eigenvalue counting function $\mathscr{N}(\tau, c_l^2-c_t^2):=\max\limits_{\substack{ k}} \{ k\big| \tau_k<\tau\}$ to eigenvalue problem (\ref{(A.9)}) or (\ref{(A.10)}). Then the two-term spectral asymptotics (\ref{(A.3)}) for linear elasticity can be equivalently written as
\begin{eqnarray} \label{(A.12)}  \mathscr{N}(\tau, c_l^2-c_t^2)= \frac{(c_l^{-2}+ c_t^{-2})S}{4\pi } \,\tau +
\frac{\beta L}{4\pi c_t}\, \tau^{\frac{1}{2}}+o(\tau^{\frac{1}{2}}), \quad \,\;\; \tau\to +\infty, \end{eqnarray}
where $S$, $L$ and $\beta$ are as in (\ref{(A.3)})--(\ref{(A.6)}).

Therefore, for each fixed $c_t^2>0$, as $c_l^2\to c_t^2$  we find that for any fixed positive integer $k\ge 1$,
\begin{eqnarray*}  \lim\limits_{c_l^2-c_t^2\to 0} \tau_k^{\mbox{Dir}}(c_l^2-c_t^2) = \tau_k^{\mbox{Dir}} (0) \quad \;\;\mbox{and}\; \;\;\;   \lim\limits_{c_l^2-c_t^2\to 0} \tau_k^{\mbox{free}}(c_l^2-c_t^2) = \tau_k^{\mbox{free}} (0),\end{eqnarray*}
 which implies
\begin{eqnarray*}&& \lim\limits_{c_l^2-c_t^2\to 0} \mathcal{Z}^{\mbox{Dir}}(t, c_l^2-c_t^2) =     \mathcal{Z}^{\mbox{Dir}}(t,0) \;\;\;\;\mbox{and} \;\;\;   \lim\limits_{c_l^2-c_t^2\to 0} \mathcal{Z}^{\mbox{free}}(t, c_l^2-c_t^2) =    \mathcal{Z}^{\mbox{free}}(t,0),\quad\;\;\\
 && \lim\limits_{c_l^2-c_t^2\to 0} \mathscr{N}^{\mbox{Dir}} (\tau, c_l^2-c_t^2)= \mathscr{N}^{\mbox{Dir}} (\tau,0)
\;\;\;\mbox{and}\;\;\;
 \lim\limits_{c_l^2-c_t^2\to 0} \mathscr{N}^{\mbox{free}} (\tau, c_l^2-c_t^2)= \mathscr{N}^{\mbox{free}} (\tau,0), \quad \;\;\;\,\end{eqnarray*}
where $$\mathcal{Z}^{\mbox{Dir}}(t,\eta):=\mbox{Tr}\; e^{-t({P}^{\mbox{Dir}}(\eta))}=\sum_{k=1}^\infty e^{-t(\tau_k^{\mbox{Dir}}(\eta))}\qquad $$  and $$\mathcal{Z}^{\mbox{free}}(t,\eta):=\mbox{Tr}\; e^{-t({P}^{\mbox{free}}(\eta))}=\sum_{k=1}^\infty e^{-t(\tau_k^{\mbox{free}}(\eta))}\qquad $$ are the {partition functions} (or the {\it traces of the heat semigroups}) with the Dirichlet and free boundary conditions, and $
{P}^{\mbox{Dir}}(\eta)$ and $
{P}^{\mbox{free}}(\eta)$ are the self-adjoint operators generated by ${P}$ with Dirichlet
and free boundary conditions, respectively.
\vskip 0.2 true cm

Further, for any fixed $c_t^2>0$, by  letting $c_l^2-c_t^2\to 0$ in (\ref{(A.9)}) and (\ref{(A.10)}), we get the following classical Laplace-type eigenvalue problems:
\begin{eqnarray}  \label{(A.13)}  \left\{ \begin{array}{ll} -c_t^2 \Delta \mathbf{v} =\tau \mathbf{v} \;\;\; &\mbox{in}\;\; M,\\
 \mathbf{v}=0 \;\;\;&\mbox{on}\;\; \partial M \end{array} \right.
\end{eqnarray}
and
\begin{eqnarray} \label{(A.14)}   \left\{ \begin{array}{ll} -c_t^2 \Delta \mathbf{v}  =\tau \mathbf{v} \;\;\; &\mbox{in}\;\; M,\\
 c_t^2 \big[\partial_{y_j} \mathsf{v}_i + \partial_{y_i}\mathsf{v}_j\big]_{n\times n}\,\boldsymbol{\nu}-c_t^2\, (\mbox{div}\; \mathbf{v})\,\boldsymbol{\nu}=0 \;\;\;&\mbox{on}\;\; \partial M. \end{array} \right.
\end{eqnarray}
We remark that in (\ref{(A.13)}) and (\ref{(A.14)}), the corresponding Rayleigh quotient  is 
\begin{eqnarray}   \label{(A. 15) }   \frac{\int_M \Big( - c_t^2 |\mbox{div}\;\mathbf{v}|^2 + \frac{c_t^2}{2}  \sum_{i,j} |\partial_{y_j} \mathsf{v}_i +\partial_{y_i} \mathsf{v}_j|^2\Big) dy}{\int_M |\mathbf{v}|^2 dy }  \end{eqnarray}
with $\mathbf{v}=(\mathsf{v}_1,\cdots, \mathsf{v}_n)\in [C_0^\infty (M)]^n$ for Dirichlet boundary condition and $\mathbf{v}\in \{\mathbf{w}\!\in\! [C^\infty (M)]^n\big| $  $\big( 2 c_t^2 (\mbox{Def}\, \mathbf{w})-c_t^2 (\mbox{div}\, \mathbf{w})\big)\,\boldsymbol{\nu}=0 \,\,\mbox{on}\,\, \partial M\}$ for free boundary condition.
It is the same as the previous argument that there exist sequences $\{(\tau_k^{\mbox{Dir}}(0), \mathbf{v}_k^{\mbox{Dir}})\}_{k=1}^\infty$ and $\{(\tau_k^{\mbox{free}}(0), \mathbf{v}_k^{\mbox{free}})\}_{k=1}^\infty$ of the eigen-pairs for the eigenvalue problems (\ref{(A.13)}) and (\ref{(A.14)}), respectively.

\vskip 0.10 true cm

(\ref{(A.13)}) is just the eigenvalue problem of the Laplace-type operator with  Dirichlet boundary condition. In fact, if $\mathbf{v}$ is an eigen-vector corresponding to $\tau$ for the eigenvalue problem (\ref{(A.13)}), then we have $\int_M \langle\mathbf{v}, -c_t^2 \Delta \mathbf{v})\rangle dy= \int_M \langle \mathbf{v}, -c_t^2 \Delta \mathbf{v} -(c_t^2-c_t^2)\, \mbox{grad}\; \mbox{div}\, \mathbf{v}\rangle dy= \int_M   \big(\! - c_t^2\, |\mbox{div}\;\mathbf{v}|^2 + \frac{c_t^2}{2}  \sum_{i,j} |\partial_{y_j} \mathsf{v}_i +\partial_{y_i} \mathsf{v}_j|^2\Big) dy$ because of $\mathbf{v}\big|_{\partial M}=0$ and $c_l^2-c_t^2=0$. On the other hand, by Green's formula we immediately have $\int_M \langle \mathbf{v}, -c_t^2 \Delta \mathbf{v}\rangle dy=c_t^2\int_M |\nabla \mathbf{v}|^2 dy$. Hence we get  $\int_M   \big(\! - c_t^2 \,|\mbox{div}\;\mathbf{v}|^2 + \frac{c_t^2}{2}  \sum_{i,j} |\partial_{y_j} \mathsf{v}_i +\partial_{y_i} \mathsf{v}_j|^2\Big) dy=c_t^2\int_M |\nabla\mathbf{v}|^2 dy$ for every eigenvector $\mathbf{v}$ of (\ref{(A.13)}).

The eigenvalue problem (\ref{(A.14)}) has a bit difference  since
 the free boundary condition in the problem (\ref{(A.14)}) (for the Laplace-type operator) is different from  the Neumann boundary condition $(\nabla u)\cdot \boldsymbol{\nu}=0$ on $\partial M$  (for the Laplace-type operator) in McKean-Singer's paper \cite{MS-67}.

For the eigenvalue problems (\ref{(A.13)}) and (\ref{(A.14)}), it is completely analogous to the discussions in \cite{MS-67} (almost verbatim, from the beginning  of \S1  to the ending of \S5 in  \cite{MS-67}, except for $e(t,x,y)$ and $e^\mp (t,x,y)$ being the matrix-valued functions and finally taking  integrals of the traces and  adding  an (auxiliary)  matrix function $-\mathbf{H}(t,x,y)$ (see below)  for the free boundary condition for the Laplacian) that
\begin{eqnarray}   \label{(A.16)}\!\!\!\!\!\!\!\!\!\!\!\!\!\!\!\!\!\!\!\!\!\!\!\!\!\sum_{k=1}^\infty e^{-t (\tau_k^{\mbox{Dir}}(0))} \,\sim  \!\!\!\!\!\!&\!\!&\!\!  \frac{n}{(4\pi c_t^2)^{n/2}\, } \,\big(\mbox{Vol}_n(M)\big)\, t^{-n/2} -  \frac{n}{4(4\pi c_t^2)^{(n-1)/2}}\,\big(\mbox{Vol}_{n-1}(\partial M)\big)\, t^{-(n-1)/2} \\
 && +O(t^{1-n/2}) \;\;\;\mbox{as}\;\; t\to 0^+,  \nonumber\end{eqnarray}
\begin{eqnarray}  \label{(A.17)} \!\! \!\!\!\!\!\!\!\! \!\!\!\!\!\!\!\!\!\!\!\!\!\!\! \sum_{k=1}^\infty e^{-t \tau_k^{\mbox{free}}(0)} \,\sim   \!\!\!\!\!\!&\!\!&\!\! \frac{n}{(4\pi c_t^2)^{n/2}\, }\,\big(\mbox{Vol}_n(M)\big)\,  t^{-n/2} +  \frac{n}{4(4\pi c_t^2)^{(n-1)/2}} \,\big(\mbox{Vol}_{n-1}(\partial M)\big)\, t^{-(n-1)/2} \\
&&+O(t^{1-n/2}) \;\;\;\mbox{as}\;\; t\to 0^+.\nonumber  \end{eqnarray}
 In fact, as in \S5 of \cite{MS-67},  let $(\Omega, g)$ be an open smooth Riemannian manifold with compact smooth boundary $\partial \Omega$,  $\,\mathscr{M}= \Omega\cup (\partial \Omega) \cup \Omega^*$ the (closed) double of $\Omega$,  and $\mathbf{Q}$ the double to $\mathscr{M}$ of the Laplace-Beltrami operator $\mathbf{I}_n \Delta_g$ 
on $\Omega$ (see p.$\,$53  in \cite{MS-67} or a detailed  explanation for $Q$ in (4.7)--(4.11) of \cite{Liu-23}), where $\mathbf{I}_n$ is the identity matrix. Define $\mathbf{Q}^-\,$ (respectively, $\mathbf{Q}^+$) to be $\mathbf{Q}\big|_{C^\infty (\bar \Omega)}$ subject to $\mathbf{u}=0\,$ (respectively, $2c_t^2  (\mbox{Def}\; \mathbf{u})^\# \, \boldsymbol{\nu} -c_t^2\, (\mbox{div}\; \mathbf{u})\,\boldsymbol{\nu}=0$) 
on $\partial \Omega$, where  $\mbox{Def} \;\! \mathbf{u}$ is the deformation tensor of the vector field $\mathbf{u}$ (see, \cite{Liu-19}),  $\#$ is the sharp operator by raising index, and $\boldsymbol{\nu}=(\nu_1,\cdots, \nu_n)$ is the unit outer normal to $\partial \Omega$.  
(Note that,  in general,  for an (elastic) vector field $\mathbf{u}$ defined in $\Omega$, the boundary traction operator can also  be equivalently written as 
$\mathcal{F} \mathbf{u} := 2c_t^2  (\mbox{Def}\; \mathbf{u})^\# \, \boldsymbol{\nu} +(c_l^2-2c_t^2)\, (\mbox{div}\; \mathbf{u})\,\boldsymbol{\nu}$ on $\partial \Omega$, i.e.,  
\begin{eqnarray} \label{23.12.12-1} \!\!\! (\mathcal{F} \mathbf{u})^k :=c_t^2 \sum_{l=1}^n\Big( \nu^{\,l} \nabla_l u^k + \nu_l \nabla^k  u^l\Big)+  (c_l^2 -2c_t^2)\, \boldsymbol{\nu}^k
\sum_{l=1}^n \nabla_l u^l \; \;\mbox{on} \;\, \partial \Omega, \,\;\; k=1,\cdots, n,\;\;\end{eqnarray}  
 where $\nabla_k u^m=u^m_{\;\;\;;\,k}\!:= \frac{\partial u^m}{\partial y_k} +  \sum_{l=1}^n \Gamma^{m}_{kl} u^l$
    is the covariant derivative of the vector field $\mathbf{u}=(u^1, \cdots, u^n)$, and $\nabla^k {u}^m\! :=u^{m;\,k}$ is the raising of  index.)    
 Let $\mathbf{E}(t, x,y)$ be the elementary solution of $\frac{\partial \mathbf{u}}{\partial t} =\mathbf{Q}\mathbf{u}$ and $\mathbf{E}^{\mp}(t,x,y)$ the elementary solution of $\frac{\partial \mathbf{u}}{\partial t}=\mathbf{Q}^{\mp} \mathbf{u}$, then 
it is easy to see $\mathbf{E} (t, x,y) = \big(E^{jk}(t,x,y)\big)_{n\times n} =  e(t,x,y))  \mathbf{I}_n$, where $e(t,x,y)$ is exactly the elementary solution of heat equation $\frac{\partial u}{\partial t} =Q u$ on $\mathscr{M}$  ($Q$ is the  
 double to $\mathscr{M}$ of the Laplace-Beltrami operator $\Delta_g$ 
on $\Omega$ (see,  p.$\,$53  in \S5 of \cite{MS-67})). Since $e(t,x,y) $ belongs to $C^\infty [ (0,\infty) \times (\mathscr{M}\setminus (\partial \Omega)) \times    (\mathscr{M}\setminus (\partial \Omega)) ] \cap C^1[ (0,\infty) \times \mathscr{M} \times \mathscr{M}\big]$ (cf.   p.$\,$53 of \cite{MS-67}). Write  $ \mathbf{E}_j(t, x, y)= (E^{j1}(t,x,y), \cdots, E^{jn} (t,x,y)$,  $\;(j=1,\cdots, n)$,   and denote   \begin{eqnarray}\label{23.12.13-2} \mathcal{F} \mathbf{E} (t,x,y) =\Big( \mathcal{F} \mathbf{E}_1 (t,x,y), 
 \cdots,  \mathcal{F} \mathbf{E}_n(t,x,y)\Big) \;\;\, \mbox{with}\;\;\, c_l^2-c_t^2=0  \end{eqnarray}  
 for all $t>0, x\in \Omega, y\in \partial \Omega$.  More precisely,  for each $j=1,\cdots, n$, the $k$-th component of $\mathcal{F} \mathbf{E}_j$ is 
 \begin{eqnarray} \label{23.12.13-6}  (\mathcal{F} \mathbf{E}_j)^k = c_t^2 \sum_{l=1}^n ( \nu^l \nabla_l E^{jk} +\nu_l \nabla^k E^{jl})
   -c_t^2 \nu^k \sum_{l=1}^n\nabla_l E^{jl}.\end{eqnarray}

Let us point out that  $\frac{\partial (\mathbf{E} (t,x,y)+ \mathbf{E}(t,x,\overset{*}{y}))}{\partial \boldsymbol{\nu}_y}\big|_{\partial \Omega} =
\frac{\partial (({e} (t,x,y)+ \mathbf{e}(t,x,\overset{*}{y}))\mathbf{I}_n)}{\partial \boldsymbol{\nu}_y}\big|_{\partial \Omega} =
0$ had been proved in \S5 of \cite{MS-67} by McKean and Singer. Now, we change $2\mathcal{F} \mathbf{E} (t,x,y)$ into another different  matrix-valued function $\boldsymbol{\Upsilon} (t,x,y)$ as following:
 changing $\frac{\partial E^{jk}(t,x,y)}{\partial \nu_y}\big|_{\partial \Omega}$ into $0$  in  the expression  of $2\mathcal{F}\mathbf{E}(t,x,y)$,  we obtain  a matrix-valued function $\boldsymbol{\Upsilon}(t,x,y)$  for $t>0$, $x\in \Omega$ and $y\in \partial \Omega$.  This implies that  $\boldsymbol{\Upsilon}(t,x,y)$ only contains the (boundary) tangent derivatives of $E^{jk}(t, x, y)$ with respect to $y\in \partial\Omega$ (without normal derivative of $E^{k}(t,x,y)$) in the local expression).  That is,  in local boundary normal coordinates (the inner normal $\boldsymbol{\nu}$ of $\partial \Omega$ is in the direction of  $x_n$-axis), \begin{eqnarray*}  \boldsymbol{\Upsilon} (t,x,y)= \Big( \boldsymbol{\Upsilon}_1 (t,x,y), \cdots,   \boldsymbol{\Upsilon}_n (t,x,y)\Big), \end{eqnarray*}
   \begin{eqnarray*}\!\!\!\!\! \!&\!\!\!&\!\!\! \boldsymbol{\Upsilon}_j (t,x,y) \\
\!\!\!\! \!\!\! &\!\!\!\!&\!\!\! \! =\!
 2c_t^2 \begin{small}   \begin{pmatrix} \! 2\frac{\partial E^{j1}}{\partial x_1}+2\!\sum\limits_{m=1}^n \!\Gamma^{1}_{1m} E^{jm}  &\!\! \cdots    & \!\!\frac{\partial E^{j1}}{\partial x_{n\!-\!1}} \!+\!\frac{\partial E^{\!j,n\!-\!1}}{\partial x_1} \!+\!\! \sum\limits_{m=1}^n\!\! \big( \Gamma^1_{\!n\!-\!1, m}\! +\!\Gamma^{n\!-\!1}_{\!1m} \big) E^{jm} 
            &   \frac{\partial E^{jn}}{\partial x_1}\!+\!\!\sum\limits_{m=1}^n\!\! \big( \Gamma^1_{\!nm} \!+\!\Gamma^n_{1m}\big) E^{jm}  
      \\
\!\! \!\cdots & \!\cdots &\!\cdots  & \cdots \\
 \!  \! \frac{\partial E^{\!j,n\!-\!1}}{\partial x_{1}}\!\!+\! \frac{\partial E^{j1}} {\partial x_{n\!-\!1}}\! +\! \!\!\sum\limits_{m=1}^n \!\! \big( \Gamma^{n\!-\!1}_{\!1m}\! +\! \Gamma^1_{\!n\!-\!1, m} \big)  E^{jm}
  &\! \cdots    &\!2 \frac{\partial E^{j,n\!-\!1}}{\partial x_{n\!-\!1}}\!+\!2\!\sum\limits_{m=1}^n \!\!\Gamma^{n\!-\!1}_{\!n\!-\!1,m} E^{jm} &   \frac{\partial E^{jn}}{\partial x_{n\!-\!1}}\!+\!\!\sum\limits_{m=1}^n\!\! \big(\Gamma^{n\!-\!1}_{\!\!nm}\! +\!\Gamma_{\!\!n\!-\!1,m}^n \big) E^{\!jm} 
   \\
  \! \frac{\partial E^{jn}}{\partial x_1} \!+\!\!\sum\limits_{m=1}^n \!\!\big( \Gamma_{\!1m}^n\! +\!\Gamma^1_{\!nm}\big) E^{jm}    &\! \cdots    & 
\frac{\partial E^{jn}}{\partial x_{n\!-\!1}}\! +\!\Gamma_{\!nm}^{n\!-\!1} E^{jm}   & 
    \;2 \sum\limits_{m=1}^n 
 \Gamma^n_{nm}  E^{jm}
 \end{pmatrix}\end{small} \begin{pmatrix} \nu_1 \\ \vdots \\ \nu_n \end{pmatrix}\\  [3mm]
\!\!\!\!&\!\!\!&\!\!\! \;- \, 2c_t^2\, \bigg( \frac{\partial E^{j1}}{\partial x_1} + \sum_{m=1}^n \Gamma^{1}_{1m} E^{jm}  + \cdots + \frac{\partial E^{j,n-1}}{\partial x_{n-1}} + \sum_{m=1}^n \Gamma^{n-1}_{n-1,m} E^{jm} +  \sum_{m=1}^n \Gamma^{n}_{nm} E^{jm} \bigg) \begin{pmatrix} \nu_1 \\ \vdots \\ \nu_n \end{pmatrix},\;\;\;  \; j=1,\cdots, n.\end{eqnarray*}
  It is easy to see that  $\boldsymbol{\Upsilon} (t,x,y)$  is a continuous (matrix-valued) function for all $t> 0$, $x\in \Omega$ and $y\in \partial \Omega$. Further, for any  fixed $x\in \Omega$ and $y\in \partial \Omega$,  since $x\ne y$ we see that 
 \begin{eqnarray} \label{23.12.18-1}  \lim\limits_{t\to 0^+} \mathbf{E} (t,x,y)=0,\end{eqnarray} 
 which implies \begin{eqnarray} \label{23.12.18-2} \lim\limits_{t \to 0^{+}} \frac{\partial \mathbf{E}(t,x,y)}{\partial T} \big|_{\partial \Omega} =0\;\;\;\mbox{for}\;\;\,   x\in \Omega, \;\;\, y\in \partial \Omega,  \end{eqnarray} 
 where $\frac{\partial}{\partial T}$ denotes the tangent derivative along the boundary $\partial \Omega$ in variable $y$. 
 From (\ref{23.12.18-1})--(\ref{23.12.18-2}) we get 
 \begin{eqnarray*} \label{23.12.18-3}  \lim\limits_{t\to 0^+} \nabla_l \mathbf{E} (t,x,y) \big|_{\partial \Omega}
 =0 \;\; \; \mbox{for all}\;\; x \in \Omega, \; y\in \partial \Omega, \;\;\, 1\le l \le n-1,\end{eqnarray*}
   so that  \begin{eqnarray} \label{23.12.18-4} \lim\limits_{t\to 0^{+}} \boldsymbol{\Upsilon}(t,x,y)=0  \;\;\,\mbox{for any } x\in \Omega, 
   \; y\in \partial \Omega,\end{eqnarray}
   where $\nabla_l \mathbf{E}= (\nabla_l \mathbf{E}_1, \cdots, \nabla_l \mathbf{E}_n)$ and $\nabla_l {E}^{jk}:= \frac{\partial 
       E^{jk} (t,x,y)}{\partial y_l}  +\sum_{m=1}^n\Gamma^k_{lm} E^{jm} (t,x,y)$.  Let $\mathbf{H} (t,x,y)$ be the solution of  
 \begin{eqnarray} \label{24.20.1} \left\{\! \begin{array}{ll}  \frac{\partial \mathbf{u} (t,x,y)}{\partial t} = \Delta_g  \mathbf{u} (t,x,y) \; \;\;\mbox{for}\;\,  t>0, \; x,y\in \Omega,\\
2  c_t^2  \big(\mbox{Def}\, \mathbf{u}(t,x,y)\big)^\#  \,\boldsymbol{\nu}  - c_t^2 \big(\mbox{div}\, (\mathbf{u}(t,x,y))\big)\, \boldsymbol{\nu} =\boldsymbol{\Upsilon}(t,x,y) \;\;\mbox{for} \,\;  t>0, \, x\in \Omega, \, y\in \partial \Omega,\\
  \mathbf{u}(0, x,y) = \mathbf{0} \,\;\;\mbox{for} \,\;  x, y \in \Omega.\end{array} \right. \;\;\;\;\end{eqnarray}
     From (\ref{23.12.18-4}), we get that the above parabolic system satisfy the compatibility condition.   
   Clearly,  the  matrix-valued function $\mathbf{H}(t,x,y)$ is  smooth in $(0, \infty) \times  \Omega \times  \Omega$ and  continuous on  $[0, \infty) \times \bar \Omega \times \bar \Omega$.  Then there exists a constant $C>0$ such that 
  \begin{eqnarray*}  | \mathbf{H} (t,x,y) |\le C \,\, \;\;\mbox{for all} \;\;   0\le t\le 1, \; x,\,y\in \bar \Omega,\end{eqnarray*}
   and hence, for dimensions $n\ge 2$,   \begin{eqnarray}\label{23.12.9-1} \int_{\Omega} \mbox{Tr}\; \mathbf{H} (t,x,x)\, dx =nC\,\mbox{vol}(\Omega) =o(t^{-\frac{n-1}{2}})\,  \;\;\mbox{as}\;\; t\to 0^+.\end{eqnarray}
(Actually, we have  $\,\lim_{t \to 0^{+}} \!\int_{\Omega} \mbox{Tr}\; \mathbf{H} (t,x,x)\, dx=0 $ by  (\ref{23.12.18-4}) and (\ref{24.20.1})).  It  can easily be verified  that 
\begin{eqnarray} 
 && \label{23.12.19-1}  \mathbf{E}^{-} (t,x,y)= \mathbf{E}(t,x,y)- \mathbf{E}(t,x,\overset{*}{y}), \\
&& \label{23.12.19-2}   \mathbf{E}^{+} (t,x,y)= \mathbf{E}(t,x,y)+ \mathbf{E}(t,x,\overset{*}{y})-  \mathbf{H}(t,x,y),
\end{eqnarray}   
where $\overset{*}{y}\in \Omega^*$ being the double of $y\in \Omega$.  
The above (\ref{23.12.19-1}) is obvious since 
$E^{jk} (t,x,y)= e(t,x,y)\, \delta_{jk}$ belongs to $C^\infty [ (0,\infty) \times (\mathscr{M}\setminus (\partial \Omega)) \times    (\mathscr{M}\setminus (\partial \Omega)) ] \cap C^1[ (0,\infty) \times \mathscr{M} \times \mathscr{M}\big]$, where $\delta_{jk}$ is the Kronicker symbol. To verify  (\ref{23.12.19-2}),  we recall that $\frac{\partial \big(\mathbf{E}(t,x,y)+\mathbf{E}(t,x,\overset{*}{y})\big)}{\partial \nu_y}\big|_{\partial \Omega}=0$
for  all $t>0$,  $x\in \Omega$ and $y\in \partial \Omega$.  
 By combining this fact and $\mathcal{F} \big(\mathbf{E} (t,x,y)+\mathbf{E}(t, x, \overset{*}{y})\big)=\boldsymbol{\Upsilon} (t,x,y)=\mathcal{F} \mathbf{H} (t,x,y)$ for $t>0$,  $x\in \Omega$ and $y\in \partial \Omega$,  we  get $\mathcal{F} \mathbf{E}^+ (t,x,y) =0$ for all $t>0$, $x\in \Omega$ and $y\in \partial \Omega$. 
  
Completely similar to the proof of (2) in \S5 of \cite{MS-67}, we can
get 
\begin{eqnarray} \label{23.12.10-1} \!\! \!\!\!\!\!\!\!  \int_\Omega \mbox{Tr}\, \Big ( \mathbf{E} (t,x,x) \mp \mathbf{E} (t, x, \overset{*}{x})  \Big) dx \!\! \!\! \!\!\! &\!\!\!&=\, \frac{n}{ (4\pi c_t^2)^{\frac{n}{2}} }\,\Big[ \mbox{the (Riemannian)  volume of} \;\, \Omega\quad \qquad \\
&& \quad \mp \frac{1}{4} \sqrt{4\pi t} \,\times \mbox{the (Riemannian) surface area of} \;\,\partial \Omega \nonumber \\
&&\quad + \frac{t}{3} \times \mbox{the curvatura integra} \,\int_{\Omega} K \nonumber  \\
 &&\quad -\frac{t}{6} \times \mbox{the integrated mean curvature} \int_{\partial \Omega}  J  + o(t^{3/2}) \Big],  \nonumber \end{eqnarray} 
 where $\int_{\partial \Omega}$ stands for the integral over $\partial \Omega$ relative to the elementary of Riemannian surface area; $K$ is the scalar curvature of $\Omega$;  and $o(t^{3/2})$ cannot be improved. 
Note that \begin{eqnarray}  \label {2.12.10-6} \sum_{k=1}^\infty e^{- t \tau_k^{\mbox{Dir}} } = \int_\Omega \mbox{Tr}\, \big( \mathbf{E}^- (t, x,x) \big) \,dx  
 = \int_\Omega  \mbox{Tr}\, \Big( \mathbf{E} (t,x,x) -\mathbf{E} (t, x, \overset{*}{x})\Big) \,dx    \end{eqnarray} 
and   \begin{eqnarray} \label{23.12.10-7}
\sum_{k=1}^\infty e^{- t \tau_k^{\mbox{free}} } 
\!\!\!\!  &&  \!\!\!= \int_\Omega \mbox{Tr}\, \big( \mathbf{E}^+ (t, x,x) \big) \,dx  \\
 && \!\!\! = \int_\Omega  \mbox{Tr} \,\Big( \mathbf{E} (t,x,x) +\mathbf{E} (t, x, \overset{*}{x})\Big) \,dx  + \int_\Omega  \mbox{Tr}\, \Big( \mathbf{H}  (t, x, x) \Big) \,dx. \nonumber \end{eqnarray} 
 From (\ref{23.12.9-1}) -- (\ref{23.12.10-7}) we  immediately obtain the previous (\ref{(A.16)}) -- (\ref{(A.17)}) for  the smooth compact Riemannian manifold $\Omega$ with smooth boundary $\partial \Omega$.  Particularly,  (\ref{(A.16)}) -- (\ref{(A.17)})  hold for the given bounded domain $M\subset  \mathbb{R}^n$ with smooth boundary $\partial M$.

   (\ref{(A.16)}) -- (\ref{(A.17)}) are very similar to the classical asymptotic formulae for the heat traces (cf. \cite{MS-67} or \cite{Gil-75}), which are correct and credible since their proofs are completely  based on McKean-Singer's celebrated method (i.e., the so-called ``method
of images'', see \cite{MS-67}). 

On the other hand, in the two-dimensional case (with condition $c_t^2>0$, $c_l^2>0$), as $c_l^2-c_t^2 \to 0^+$ (without loss of generality, $c_t^2>0$ can be fixed) one has  $\alpha= c_t^2c_l^{-2} \to 1^-$ (As pointed out earlier, one can take $\alpha=1$ in \cite{SaVa-97} because the Lam\'{e} operator (with both boundary conditions) is  strongly elliptic for $\alpha=1$ and  uniformly strongly elliptic for all $\alpha \in [\alpha_0, \alpha_1]$, where $\alpha_0$ and $\alpha_1$ are  arbitrary positive numbers with $\;0<\alpha_0<\alpha <\alpha_1<+\infty$). For $\alpha=1$,  we  get that the algebraic equation (\ref{(A.6)}) (simply denoted by $R_\alpha (\gamma)=0$) has the following roots:
\begin{eqnarray}   \label{(A.18)}  \left.\begin{array}{ll}\gamma_1=0, \;\; \;\;\gamma_2= \sqrt{4+2\sqrt{2}},\; \;\; \;\gamma_3= - \sqrt{4+2\sqrt{2}},\\
[1.5mm] \gamma_4= \sqrt{4-2\sqrt{2}},\;\;\;\; \gamma_5= -\sqrt{4-2\sqrt{2}}\end{array}\right.,\end{eqnarray}
where the root $\gamma_1=0$ has multiplicity $2$.
 It follows from (\ref{(A.4)}) and (\ref{(A.5)})  that   $\beta^{\mbox{Dir}}=-2$ for the Dirichlet boundary condition as $\alpha\to 1^-$, and
\begin{eqnarray}   \label{(A.19)}   \left.\begin{array}{ll} & \beta_1^{\mbox{free}}=+\infty, \;\;\;\; \beta_2^{\mbox{free}}=-2+ \frac{4}{\sqrt{4+2\sqrt{2}}},\;\;\; \;
\beta_3^{\mbox{free}}= -2 -  \frac{4}{\sqrt{4+2\sqrt{2}}}, \\
 & \beta_4^{\mbox{free}}=-2+ \frac{4}{\sqrt{4-2\sqrt{2}}},\;\;\; \; \beta_5^{\mbox{free}}=-2- \frac{4}{\sqrt{4-2\sqrt{2}}}\end{array} \right. \end{eqnarray}
for the free boundary condition as $\alpha\to 1^-$.

\vskip 0.15 true cm

It is well-known that  the {\it partition functions}  $\mathcal{Z}^{\mbox{Dir}}(t):=\sum_{k=1}^\infty e^{-t\tau_k^{{}_{\mbox{Dir}}}}$  and  $\mathcal{Z}^{\mbox{free}}(t):=\sum_{k=1}^\infty e^{-t\tau_k^{{}_{\mbox{free}}}}$ are  just
 the Riemann-Stieltjes integrals of $e^{-t\lambda}$ with respect to the counting functions $\mathscr{N}^{\mbox{Dir}}(\lambda)$ and $\mathscr{N}^{\mbox{free}}(\lambda)$, respectively.
 That is, \begin{eqnarray} \label{(A.20)}  \mathcal{Z}^{\mbox{Dir}}(t) = \int_{-\infty}^{+\infty} e^{-t\lambda } d\mathscr{N}^{\mbox{Dir}} (\lambda), \;\;\;\;\; \mathcal{Z}^{\mbox{free}}(t) = \int_{-\infty}^{+\infty} e^{-t\lambda } d\mathscr{N}^{\mbox{free}} (\lambda). \end{eqnarray}
Thus, if the following two-term spectral asymptotics hold
 \begin{eqnarray} \label{(A.21)} \!\!\!\!\!\!\!\!\!\!\!\!\mathscr{N}^{\mbox{Dir}} (\lambda) = a \big(\mbox{Vol}_n (M)\big)\, \lambda^{n/2} + b^{\mbox{Dir}}\big( \mbox{Vol}_{n-1} (\partial M)\big) \lambda^{(n-1)/2}+ o(\lambda^{(n-1)/2}) \;\; \mbox{as}\;\; \lambda\to +\infty,\end{eqnarray}
\begin{eqnarray}  \label{(A.22)} \!\!\!\!\!\!\!\!\!\!  \mathscr{N}^{\mbox{free}} (\lambda) = a \,\big(\mbox{Vol}_n (M)\big) \lambda^{n/2} + b^{\mbox{free}}\,\big( \mbox{Vol}_{n-1} (\partial M)\big) \lambda^{(n-1)/2}+ o(\lambda^{(n-1)/2}) \;\; \mbox{as}\;\; \lambda\to +\infty, \end{eqnarray}
then, one immediately find by using (\ref{(A.20)}), (\ref{(A.21)}) and (\ref{(A.22)}) that
\begin{eqnarray}\label{(A.23)} \!\!\!\!\!\!\!\!\!\! \mathcal{Z}^{\mbox{Dir}}(t) = c\, \mbox{Vol}_n( M)\, t^{-n/2} + d^{\mbox{Dir}}\,\big( \mbox{Vol}_{n-1} (\partial M) \big)\,t^{-(n-1)/2} + o(t^{-(n-1)/2})\;\; \mbox{as} \;\; t\to 0^+,\end{eqnarray}
\begin{eqnarray}\label{(A.24)}\!\!\!\!\!\!\!\!\!\!  \mathcal{Z}^{\mbox{free}}(t) = c\, \mbox{Vol}_n( M)\, t^{-n/2} + d^{\mbox{free}}\,\big( \mbox{Vol}_{n-1} (\partial M)\big) \,t^{-(n-1)/2} + o(t^{-(n-1)/2})\;\; \mbox{as} \;\; t\to 0^+,\end{eqnarray}
   where  \begin{eqnarray} \label{(A.25)} c = \Gamma \Big( 1+\frac{n}{2}\Big)\, a, \, \; \,\;  d^{\mbox{Dir}} = \Gamma\Big(1+\frac{n-1}{2}\Big) b^{\mbox{Dir}}, \;\;\;\, d^{\mbox{free}} = \Gamma\Big(1+\frac{n-1}{2}\Big) b^{\mbox{free}}. \end{eqnarray}

\vskip 0.2 true cm

Finally, for $n=2$ (i.e., the two-dimensional case), suppose by contradiction that the result in \cite{SaVa-97} (i.e., the above spectral asymptotic formula (\ref{(A.3)}) along with (\ref{(A.5)}))  is correct for the free boundary condition. Then, (\ref{(A.12)})  (along with (\ref{(A.5)})) holds when $c_l^2-c_t^2=0$ with $c_t^2>0$. From this, (\ref{(A.20)})--(\ref{(A.25)}) and $\alpha=c_t^2c_l^{-2}=1$,  one immediately obtains the following heat trace asymptotic expansion:
  \begin{eqnarray}\label{(A.26)} \qquad \qquad \quad && \sum_{k=1}^\infty e^{-t \tau_k^{\mbox{free}}(0)} \sim  \frac{2S}{(4\pi c_t^2)\, } \,t^{-1} \!+\! \frac{L}{ 4\sqrt{\pi}\, c_t}  \Big[1 +\big(\frac{4}{\gamma_j}-3\big)\Big]\, t^{-1/2}\qquad \qquad \qquad \qquad \\
 && \qquad \qquad \quad\qquad \; \;\; +o(t^{-1/2}) \;\;\mbox{as}\;\; t\to 0^+,  \;\;\;\;\, j=1,2,3,4,5,\nonumber \end{eqnarray}
 which is different from generalized McKean-Singer's  result (\ref{(A.17)}) for the free boundary condition when $n=2$ (here, i.e., in the result of \cite{SaVa-97}, a superfluous constant $ \frac{L}{4\sqrt{\pi}\, c_t} \Big[\frac{4}{\gamma_j} -3\Big]$ appears), where $\gamma_j\,$  ($j=1,2,3,4,5$) are given by (\ref{(A.18)}). The reason is that one always has  $\frac{4}{\gamma_j} -3\ne 0$ for each $j=1,2,3,4,5$ because of (\ref{(A.18)}).
{This  implies that the elastic spectral asymptotic formula (\ref{(A.3)}) (along with  (\ref{(A.5)})) for free boundary condition is wrong. In other words, the conclusion on p.$\,$237 in \cite{SaVa-97}  is wrong.}

\vskip 0.2 true cm

With a similar method, we can also show that  the two-term spectral asymptotic formula for three-dimensional elasticity  with free boundary condition (see, p.$\,$237--238  in \cite{SaVa-97}) is also wrong.

\vskip 0.52 true cm

 \noindent{\bf Remark 2.3.} \  (i) \ These fundamental errors for the elastic spectral asymptotics formulae stem from the so-called ``algorithm'' in \cite{SaVa-97}, which, throughout the whole book \cite{SaVa-97}, is essentially  wrong. It is impossible to correct these kinds of fundamental errors at all. The book \cite{SaVa-97} has misled  a large number of readers for twenty-six years.

(ii) \   {The same fundamental  errors also appeared in a series of papers \cite{CaFrLeVa-23}, \cite{Va-84}, \cite{Va-86}}. For more detailed arguments, we refer the reader to \cite{Liu-23}.

(iii) \  {Clearly, most conclusions in the book \cite{SaVa-97} are wrong because they  are based on an erroneous ``algorithm'' method}.

\vskip 0.69 true cm  
\section{Heat trace asymptotic expansion for the Lam\'{e} operator on a Riemannian manifold}
\vskip 0.19 true cm  

Let $(\Omega,g)$ be a compact, smooth, $n$-dimensional Riemannian manifold with smooth boundary $\partial \Omega$. Let $P_g$ be the Navier--Lam\'{e} operator (see \cite{Liu-19} and \cite{Liu-21}):
\begin{eqnarray} \label {1-1} P_g\mathbf{u}:=c_t^2 \nabla^* \nabla \mathbf{u} -(c_l^2-c_t^2) \,\mbox{grad}\; \mbox{div}\, \mathbf{u} -c_t^2 \, \mbox{Ric} (\mathbf{u}), \;  \;\; \mathbf{u}=(u^1, \cdots, u^n),\end{eqnarray}
where  the Lam\'{e} coefficients (constants) $c_t^2$ and $c_l^2-2c_t^2$ still satisfy $c_t^2>0$ and $c_l^2>  0$,  $\nabla^* \nabla $ is the Bochner Laplacian (see (2.11) of \cite{Liu-21}), $\mbox{div}$ and $\mbox{grad}$ are the usual divergence and gradient operators, and  \begin{eqnarray} \label{18/12/22} \mbox{Ric} (\mathbf{u})= \big(\sum\limits_{l=1}^nR^{\,1}_{l} u^l,  \sum\limits_{l=1}^n R^{\,\,2}_{\,l} u^l, \cdots, \sum\limits_{l=1}^n R^{\,\,n}_{\,l} u^l\big)\end{eqnarray} denotes the action of Ricci tensor $\mbox{R}_l^{\;j}:=\sum_{k=1}^n R^{k\,\,j}_{\,lk}$ on $\mathbf{u}$.
We denote by $P_g^-$ and $P_g^+$ the Navier--Lam\'{e} operators with zero Dirichlet and zero Neumann boundary conditions, respectively. The Dirichlet boundary condition is $\mathbf{u}\big|_{\partial \Omega}$, and the  free  boundary condition is \begin{eqnarray*}2c_t^2  (\mbox{Def}\, \mathbf{u})^\# \boldsymbol{\nu} + (c_l^2-2c_t^2)  (\mbox{div}\, \mathbf{u} )\boldsymbol{\nu}\;\;\mbox{on}\;\; \partial \Omega,\end{eqnarray*} where $\mbox{Def}\, \mathbf{u}= \frac{1}{2} (\nabla \mathbf{u} +(\nabla \mathbf{u})^T)$, $\,(\nabla \mathbf{u})^T$ is the transpose of $\nabla \mathbf{u}$. $ \;\#$ is the sharp operator (for a tensor) by raising index, and $\boldsymbol{\nu}$ is the unit outer normal to $\partial \Omega$.
Since $P_g^-$ (respectively, $P_g^+$) is an unbounded, self-adjoint and positive (respectively,  nonnegative) operator in $[H^1_0(\Omega)]^n$ (respectively, $[H^1(\Omega)]^n$) with discrete spectrum $0< \tau_1^- < \tau_2^- \le \cdots \le \tau_k^- \le \cdots \to +\infty$ (respectively, $0\le \tau_1^+ < \tau_2^+ \le \cdots \le \tau_k^+ \le \cdots \to +\infty$), one has
\begin{eqnarray} \label{1-4} P_g^\mp {\mathbf{u}}_k^\mp =\tau_k^\mp {\mathbf{u}}_k^\mp,\end{eqnarray} where ${\mathbf{u}}_k^-\in  [H^1_0(\Omega)]^n$ (respectively, ${\mathbf{u}}_k^+\in  [H^1(\Omega)]^n$) is the eigenvector corresponding to elastic eigenvalue $\tau_k^{-}$ (respectively, $\tau_k^{+}$).

\vskip 0.10 true cm

Combining the previous discussions  for the range of the Lam\'{e} coefficients (i.e.,  the strongly ellipticity conditions $c_t^2>0$ and $c_l^2>0$), the method of the heat trace and  ``method of image'' , similar to \cite{Liu-21}  (also see \cite{Liu-23})    we can obtain the following result 
(i.e.,  two-term asymptotic expansions for the traces of (elastic) heat semigroups):

\vskip 0.22 true cm

\noindent{\bf Theorem 3.1.} \ {\it
Let $(\Omega,g)$ be an $n$-dimensional compact smooth  Riemannian manifold with smooth boundary $\partial \Omega$.  Assume that  the Lam\'{e} coefficients $c_t^2$ and $c_l^2-2 c_t^2$ satisfy $c_t^2>0$ and $c_l^2> 0$.  Let $0< \tau_1^-< \tau_2^- \le \tau^-_3\le \cdots \le \tau_k^- \le \cdots$ (respectively, $0\le \tau_1^+ < \tau_2^+ \le \tau_3^+ \le \cdots \le \tau_k^+ \le \cdots $) be the eigenvalues of the Navier--Lam\'{e} operator $P_g^-$ (respectively, $P_g^+$) with respect to the zero Dirichlet (respectively, zero Neumann) boundary condition. Then
\begin{eqnarray}  \mathcal{Z}^{\mp}(t) \! \!\!\!\! &&\!\!\!=  \sum_{k=1}^\infty e^{-t \tau_k^{\mp}} =\bigg[ \frac{n-1}{(4\pi c_t^2 t)^{n/2}}
 + \frac{1}{(4\pi c_l^2 t)^{n/2}}\bigg] {\mbox{Vol}_n}(\Omega) \nonumber\\
\!\!\!\!\!&& \;\label{23.1-7} \, \mp \frac{1}{4} \bigg[  \frac{n-1}{(4\pi c_t^2 t)^{(n-1)/2}}
 +  \frac{1}{(4\pi c_l^2 t)^{(n-1)/2}}\bigg]{\mbox{Vol}_{n-1}}(\partial\Omega)+O(t^{{1-n}/2})\quad \mbox{as}\;\; t\to 0^+.\quad \quad \;\;\end{eqnarray}
 Here ${\mbox{Vol}}_{n}(\Omega)$ denotes the $n$-dimensional volume of $\, \Omega$,  ${\mbox{Vol}}_{n-1} (\partial\Omega)$ denotes the $(n-1)$-dimensional volume of $\partial \Omega$,  
  $c_t$ is  the velocity of transverse wave which is a constant (here the subscript $t$ of $c_t^2$ is the first letter of the word ``transverse'', it is not the time variable). }

\vskip 0.36  true cm

 {Hence, in the $n$-dimensional case ($n\ge 2$), our result  (\ref{23.1-7})  (see also, \cite{Liu-21}) is compatible with the classical McKean-Singer's result (\ref{(A.16)})--(\ref{(A.17)}) as $c_l^2-c_t^2\to 0$ (with $c_t^2>0$ is fixed).  However, the corresponding results in \cite{SaVa-97} and \cite{CaFrLeVa-23}  contradict with the classical McKean-Singer's result as $c_l^2-c_t^2\to 0$. That is, the main result of  \cite{CaFrLeVa-23} and the most conclusions in \cite{SaVa-97} are wrong. This implies the strategy adopted in \cite{SaVa-97} and \cite{CaFrLeVa-23} (originated from \cite{Va-84} and \cite{Va-86}) is fundamentally flawed. From (\ref{23.1-7}), we also see that 
 \begin{eqnarray*} d^{\mbox{Dir}} + d^{\mbox{free}}\!\!\!\!\!&&\!\!\! = -  \frac{1}{4} \bigg[  \frac{n-1}{(4\pi c_t^2)^{(n-1)/2}}
 +  \frac{1}{(4\pi c_l^2 )^{(n-1)/2}}\bigg]{\mbox{Vol}_{n-1}}(\partial\Omega)
 \\
&&\;\,  +  \frac{1}{4} \bigg[  \frac{n-1}{(4\pi c_t^2)^{(n-1)/2}}
  +  \frac{1}{(4\pi c_l^2 )^{(n-1)/2}}\bigg]{\mbox{Vol}_{n-1}}(\partial\Omega) =0,  \end{eqnarray*}
  where $d^{\mbox{Dir}}$ and $d^{\mbox{free}}$ are  as in (\ref{(A.23)}) and (\ref{(A.24)}). 
 It is well known that an asymptotic expansion of the eigenvalue counting function can infer the corresponding asymptotic expansion of the  trace  of elastic heat semigroup, and the previous (\ref{(A.25)}) must hold. Therefore,  in the asymptotic expansions of the eigenvalue counting functions with the Dirichlet and free boundary conditions, the sum of the second coefficients for two kinds of boundary conditions must vanish, i.e., $b^{\mbox{Dir}}+b^{\mbox{free}}=0$, where $b^{\mbox{Dir}}$ and $b^{\mbox{free}}$ are as in (\ref{(A.21)}) and (\ref{(A.22)}). Clearly, the sum of (\ref{(A.4)}) and (\ref{(A.5)}) (i.e., the sum of the corresponding results in \cite{SaVa-97}) does not vanish.

\vskip 0.52 true cm

 \noindent{\bf Remark 3.2.}  \  (i) \  In \cite{Liu-21}, we established an algorithm  method of  pseudodifferential operators, by which all coefficients $a_k^{\mp}$ can be explicitly calculated in the asymptotic expansions $\mathcal{Z}^{\mp}(t) \sim \sum_{k=0}^\infty a^{\mp}_k t^{-(n-k)/2}$ as $t\to 0^+$ .   
  
 (ii) \ For an $n$-dimensional Riemannian manifold $(\Omega, g)$, in two-term asymptotic expansion of the trace $\mathcal{Z}(t)$ of the elastic heat semigroups,  it is enough only if this manifold $(\Omega, g)$ is smooth compact with smooth boundary. The reason is that the initial boundary value problem of an elastic parabolic system $(\frac{\partial    }{\partial t} -P_g)\mathbf{u}(t,x)=0$ (or equivalently, the elastic heat semigroup) has quite good regularity.  
 In general, the validity of two-term asymptotic expansion for the (elastic)  eigenvalue counting function $\mathscr{N}(\lambda)$ is still an open
question (as it is for the scalar Dirichlet or Neumann Laplacian). In order to get a two-term  asymptotic expansion  for $\mathscr{N}(\lambda)$,  besides $(\Omega,g)$ being smooth compact with smooth boundary, an additional assumption (i.e., the corresponding billiards is neither deadend nor absolutely periodic) should be satisfied (see, \cite{Liu-23} for the detailed remark) .

\vskip 1.06 true cm

\section*{Acknowledgments}

 This research was supported by the NNSF of China (12271031)
 and  NNSF of China (11671033/A010802).

\vskip 0.89 true cm

\end{document}